\documentclass[12pt]{amsart}
\usepackage{amssymb}
\usepackage[english]{babel}
\input xy
\usepackage[all]{xy}

\usepackage{amsmath}

\usepackage[latin1]{inputenc}%
\usepackage{a4}

\newlength{\itemlaenge}

\newtheoremstyle{mytheorem}
  {}
  {}
  {\slshape}
  {}
  {\scshape}
  {.}
  { }
  {}

\newtheoremstyle{mydefinition}
  {}
  {}
  {\upshape}
  {}
  {\scshape}
  {.}
  { }
  {}

\theoremstyle{mytheorem}
\newtheorem{lemma}{Lemma}[section]

\newtheorem*{prop*}{Proposition}

\newtheorem{cor}[lemma]{Corollary}

\newtheorem{thm}[lemma]{Theorem}

\newtheorem*{thm*}{Theorem}
\theoremstyle{mydefinition}
\newtheorem{rem}[lemma]{Remark}
\newtheorem*{rem*}{Remark}

\newtheorem*{notation*}{Notation}

\newtheorem*{warning*}{Warning}

\newtheorem{defi}[lemma]{Definition}
\newtheorem*{defi*}{Definition}

\numberwithin{equation}{section}

\newcommand{\bqn}{\begin{eqnarray*}}
\newcommand{\eqn}{\end{eqnarray*}}
\newcommand{\bq}{\begin{eqnarray}}
\newcommand{\eq}{\end{eqnarray}}
\newcommand{\ba}{\begin{aligned}}
\newcommand{\ea}{\end{aligned}}
\newcommand{\be}{\begin{enumerate}}
\newcommand{\ee}{\end{enumerate}}

\newcommand{\thismonth}{\ifcase\month 
  \or January\or February\or March\or April\or May\or June%
  \or July\or August\or September\or October\or November%
  \or December\fi}

\newcommand{\Sp}{\operatorname{Sp}}

\newcommand{\PSL}{\operatorname{PSL}}

\newcommand{\PU}{\operatorname{PU}}

\newcommand{\CC}{{\mathbb C}}
\newcommand{\DD}{{\mathbb D}}

\newcommand{\RR}{{\mathbb R}}

\newcommand{\ZZ}{{\mathbb Z}}

\newcommand{\Xx}{{\mathcal X}}

\newcommand{\rad}{{\mathrm{Rad}}}

\def\h{{\rm H}}

\def\one{\mathbf{1\kern-1.6mm 1}}

\def\h2{{\operatorname{H_2}}}
\def\h1{{\operatorname{H_1}}}

\def\PSL{\operatorname{PSL}}

\def\Sp{\operatorname{Sp}}

\def\cs{{\check S}}

\def\kgb{\kappa_G^{\rm b}}

\def\to{\rightarrow}

\def\h{{\rm H}}

\renewcommand{\phi}{\varphi}

\def\No{N\raise4pt\hbox{\tiny o}\kern+.2em}
\def\no{n\raise4pt\hbox{\tiny o}\kern+.2em}

\renewcommand{\hom}{\textup{Hom}}

\begin{document}
\bibliographystyle{plain}
\title{Weakly maximal representations of surface groups}
\author[G.~Ben Simon]{Gabi Ben Simon}
\email{gabi.ben-simon@math.ethz.ch}
\address{Department Mathematik, ETH Zentrum, R\"amistrasse 101, CH-8092
Z\"urich, Switzerland}
\author[M.~Burger]{Marc Burger}
\email{burger@math.ethz.ch}
\address{Department Mathematik, ETH Zentrum, R\"amistrasse 101, CH-8092
Z\"urich, Switzerland}
\author[T.~Hartnick]{Tobias Hartnick}
\email{tobias.hartnick@gmail.com}
\address{Mathematics Department, Technion - Israel Institute of Technology,
Haifa, 32000, Israel}
\author[A.~Iozzi]{Alessandra Iozzi}
\email{iozzi@math.ethz.ch}
\address{Department Mathematik, ETH Zentrum, R\"amistrasse 101,
  CH-8092 Z\"urich, Switzerland}
\author[A.~Wienhard]{Anna Wienhard}
\email{wienhard@math.princeton.edu}
\address{Department of Mathematics, Princeton University, Fine Hall, Washington
Road, Princeton, NJ 08540, USA}
\thanks{T.~H. was supported in part at the Technion by a fellowship of the Israel Council of Higher Education.
A.~I. was partially supported by the Swiss National Science Foundation
project 2000021-127016/2; A.~W. was partially supported by the National Science
Foundation under agreement No.s DMS-0803216, DMS-1107367 and DMS-1065919.}


\date{\today}

\begin{abstract} 
We introduce and study a new class of representations of surface groups into Lie
groups of Hermitian type, 
called {\em weakly maximal} representations.  They are defined in terms of
invariants in bounded cohomology and extend
considerably the scope of maximal representations studied in
\cite{Burger_Iozzi_Wienhard_toledo, Burger_Iozzi_Wienhard_anosov,
Burger_Iozzi_Labourie_Wienhard, Burger_Iozzi_Wienhard_survey, Wienhard_mapping,
Hartnick_Strubel, Guichard_Wienhard_invariants, Gothen, Bradlow_GarciaPrada_Gothen_survey, Bradlow_GarciaPrada_Gothen,
Bradlow_GarciaPrada_Gothen_sp4, GarciaPrada_Gothen_Mundet}.  We prove that weakly maximal representations
are discrete and injective and describe the structure of the Zariski closure of
the image.  An interesting feature of these representations
is that they admit an elementary topological characterization in terms of
bi-invariant orderings.  In particular if the target group
is Hermitian of tube type, the ordering can be described in terms of the causal
structure on the Shilov boundary.  
\end{abstract}
\maketitle
\centerline{--- Research Announcement ---}

\section{Introduction}
This research announcement presents results obtained by the authors during the
last two years concerning the class of so-called {weakly maximal
representations}
of surface groups into a
Hermitian Lie group. A more detailed version of this note with full proofs is
currently under preparation \cite{BenSimon_Burger_Hartnick_Iozzi_Wienhard}.\\ 

Given a compact oriented surface $\Sigma$ of negative Euler characteristic,
possibly with boundary, a
general theme is to study the space of representations
$\hom(\pi_1 (\Sigma), G)$ of the fundamental group of $\Sigma$ into a semisimple
Lie group $G$, and in
particular to distinguish subsets of geometric significance, such as holonomy
representations of
geometric structures. Classical examples include the set of Fuchsian
representations in  $\hom(\pi_1
(\Sigma), \PSL(2,\RR))$ or the set of quasi-Fuchsian representations in
$\hom(\pi_1 (\Sigma), \PSL(2,\CC))$, where the target group is of real rank one.
In recent years these studies have been extended to the case where $G$ is of
higher rank. Prominent examples of geometrically significant subsets of
representation varieties for higher rank targets include Hitchin
components \cite{Hitchin, Goldman_Choi, Labourie_anosov}, positive
representations \cite{Fock_Goncharov}, maximal representations
\cite{Burger_Iozzi_Wienhard_toledo, Burger_Iozzi_Wienhard_anosov,
Burger_Iozzi_Labourie_Wienhard, Burger_Iozzi_Wienhard_survey, Wienhard_mapping,
Hartnick_Strubel, Gothen, Bradlow_GarciaPrada_Gothen_survey, Bradlow_GarciaPrada_Gothen,
Bradlow_GarciaPrada_Gothen_sp4, GarciaPrada_Gothen_Mundet} and Anosov representations
\cite{Guichard_Wienhard_anosov, Labourie_anosov, Guichard_Wienhard_invariants}. Even though these subsets
exhibit several common properties, which has lead to
summarizing their study under the terminus \emph{higher Teichm\"uller theory},
they
are defined and investigated by very different methods.\\

The present article is concerned with an extension of higher Teichm\"uller
theory in the Hermitian context. Recall that a semisimple Lie group $G$ is
called \emph{Hermitian} if the associated symmetric space $\Xx$ admits a
$G$-invariant K\"ahler form $\omega_{\Xx}$. This K\"ahler form can be used
to define a continuous function 
$T: \hom(\pi_1 (\Sigma), G) \to \RR$ on the representation variety; the
invariant $T(\rho)$ is called the \emph{Toledo number} of the
representation $\rho$. If the surface $\Sigma$ is closed then $T$ is defined
by the formula
\[T(\rho) := \frac{1}{2\pi} \cdot \int_{\Sigma} p_*f^*\omega_{\Xx},\]
where $f: \widetilde{\Sigma} \to \Xx$ is an arbitrary $\rho$-equivariant map and
$p: \widetilde{\Sigma} \to \Sigma$ is the universal covering projection. For
$\Sigma$ with boundary a modification of this definition has been provided in
\cite{Burger_Iozzi_Wienhard_toledo}.
In any case, the Toledo number is subject to a \emph{Milnor-Wood type
inequality} of the form
\begin{eqnarray}\label{MW}
|T(\rho)|\leq ||\kgb|| \cdot |\chi(\Sigma)|,
\end{eqnarray}
where $\kgb \in H^2_{cb}(G; \RR)$ denotes the bounded K\"ahler class of $G$,
i.e. the class corresponding to $\omega_\Xx$ under the isomorphisms $H^2_{cb}(G;
\RR) \cong H^2_{c}(G; \RR) \cong \Omega^2(\Xx)^G$, and $||\cdot||$ denotes the
seminorm in continuous bounded cohomology (see
\cite{Burger_Iozzi_Wienhard_toledo}). The class of representations $\rho$ with
maximal Toledo invariant $T(\rho) = ||\kgb|| \cdot |\chi(\Sigma)|$, or
\emph{maximal representations} for short, has been the main object of study in
higher Teichm\"uller
theory with Hermitian target groups
(\cite{Burger_Iozzi_Wienhard_toledo, Burger_Iozzi_Wienhard_anosov,
Burger_Iozzi_Labourie_Wienhard, Burger_Iozzi_Wienhard_survey, Wienhard_mapping,
Hartnick_Strubel, Guichard_Wienhard_invariants, Gothen, Bradlow_GarciaPrada_Gothen_survey, Bradlow_GarciaPrada_Gothen,
Bradlow_GarciaPrada_Gothen_sp4, GarciaPrada_Gothen_Mundet}). Here we propose a generalization of maximal
representations, which preserves many of their key properties. Our starting
point is the observation that the inequality \eqref{MW}
can be refined into the chain of inequalities
\begin{eqnarray*}
 |T(\rho)|\leq ||\rho^*\kgb|| \cdot |\chi(\Sigma)| \leq ||\kgb|| \cdot
|\chi(\Sigma)|.
\end{eqnarray*}
In particular, a representation is maximal iff it satisfies both $||\rho^*\kgb||
= ||\kgb||$
and $T(\rho) = ||\rho^*\kgb|| \cdot |\chi(\Sigma)|$.
Representations satisfying $||\rho^*\kgb|| = ||\kgb||$ are called \emph{tight};
these have been investigated in much greater generality in
\cite{Burger_Iozzi_Wienhard_tight}. Here we are interested in representations
satisfying the complementary property (see \cite{Wienhard_thesis}):
\begin{defi}
A representation $\rho: \pi_1(\Sigma) \to G$ is \emph{weakly maximal} if it
satisfies 
\begin{eqnarray}\label{WMDefIneq}
 T(\rho)= ||\rho^*\kgb|| \cdot |\chi(\Sigma)|.
\end{eqnarray}
\end{defi}
By definition a representation is maximal iff it is weakly maximal and tight.
In the following three sections we will present results concerning
\begin{itemize}
 \item structure theorems describing the range, kernel and Zariski closure of
weakly maximal representations with nonzero Toledo invariant;
 \item a geometric interpretation of weakly maximal representations of nonzero
Toledo invariant in terms of causal structures on Shilov boundaries;
 \item the relation between the space of weakly maximal representations and
other prominent subsets of the representation variety.
\end{itemize}
These results demonstrate that weakly maximal representations form a very broad
and geometrically significant class of representations which still share many
desirable structural properties with maximal representations.

\section{Structure theorems for weakly maximal representations}
Various general structure theorems for maximal representations have been
established by three of the present authors in
\cite{Burger_Iozzi_Wienhard_toledo}. We show that an essential part of these
structure theorems can be established for
weakly maximal representations. For the proofs it is important to understand
the range of the Toledo number when restricted to weakly maximal
representations. Let us assume that the K\"ahler form $\omega_{\Xx}$ has been
normalized to have minimal holomorphic curvature $-1$. With this normalization
we then have 
$||\kgb|| |\chi(\Sigma)| \in \ZZ$, hence maximal representations have integral
Toledo number. Moreover, if $\Sigma$ is closed, then 
$T(\rho) \in e_G^{-1} \cdot \ZZ$ for some integer $e_G$ depending only on $G$.
In particular, the range of $T$ is finite. In strong contrast, if
the surface is admitted to have a boundary, then $T(\rho)$ can take arbitrary
values inside the closed interval $[-||\kgb|| |\chi(\Sigma)| ,
||\kgb|| |\chi(\Sigma)|]$. It is therefore significant that for weakly maximal
representations we can prove:
\begin{thm}\label{thm_intro:integer}
There is natural number $n_G$ depending only on $G$, such that for every weakly
maximal representation $\rho: \pi_1(\Sigma) \to G$ we have $n_G T(\rho) \in
\ZZ$. In particular, $T$ takes only finitely many values on weakly-maximal
representations.
\end{thm}
Theorem \ref{thm_intro:integer} plays a crucial part in establishing the
following basic properties of weakly maximal representations of nonzero Toledo
invariant:
\begin{thm}\label{thm_intro: discrete_faithful}
Let $\rho: \pi_1(\Sigma) \to G$ be a weakly maximal representation and $T(\rho)
\neq 0$. Then $\rho$ is faithful with discrete image.
\end{thm}
An important step in the proof of Theorem \ref{thm_intro: discrete_faithful} is
the realization that a representation $\rho$ is weakly maximal iff there exists
$\lambda \geq 0$ such that 
\begin{eqnarray}\rho^*\kgb = \lambda \cdot \kappa_\Sigma^b,\end{eqnarray}
where $\kappa_\Sigma^b \in H^2_b(\Gamma)$ is the bounded fundamental class of
the surface $\Sigma$ as introduced in \cite{Burger_Iozzi_Wienhard_toledo}. By
Theorem \ref{thm_intro:integer} the constant $\lambda$ has in fact to be
\emph{rational}. This provides severe restrictions on the kernel and range of
$\rho$.\\

Both Theorem \ref{thm_intro:integer} and Theorem \ref{thm_intro:
discrete_faithful} depend on understanding the Zariski closure of a weakly
maximal representation. Unlike for maximal representations, the Zariski closure
of a weakly maximal representation need not be reductive. To overcome this
difficulty, we argue as follows. We first show that for a closed subgroup
$L < G$ there exists a unique maximal
normal subgroups of $L$  on which $\kgb|_L$ vanishes. This subgroup is called
the {\em K\"ahler radical} $\rad_{\kgb}(L)$of $L$, and the quotient
$L/\rad_{\kgb}(L)$ is automatically semisimple. While the Zariski closure of a
weakly maximal representation can be fairly complicated, we have a rather good
control over its quotient by its K\"ahler radical, provided the Toledo number is
nonzero: 
\begin{thm}\label{ZariskiClosure}
Let $\rho: \pi_1(\Sigma) \to G$ be weakly maximal representation with
$T(\rho)\neq 0$. Let $L<G$ be the Zariski closure of the image of $\rho$ and set
$H= L/ \rad_{\kgb}(L)$. 
Then 
\begin{enumerate}
\item  $H$ is a semisimple Lie group of Hermitian type; all almost simple factor
of $H$ are of tube type. 
\item The composition $\pi_1(\Sigma) \to L \to H$ is faithful with discrete
image. 
\end{enumerate}
\end{thm}
\begin{rem}
In the above theorems it is essential that the Toledo number is nonzero. However
the class of weakly maximal representations with $T(\rho) = 0$ is also of
interest. It includes in particular the set of representations where
$\rho^*(\kgb) = 0$. In the case when $G = \PU(1,n)$ such representations have
been studied in \cite{Burger_Iozzi_totallyreal}.
\end{rem}

\section{Geometric description of weakly maximal representations}
It turns out that techniques from \cite{BenSimon_Hartnick_Comm} can be used to
provide
a geometric characterization of weakly maximal representations with nonzero
Toledo invariant in terms of bi-invariant orders. 
To simplify the formulation we will only spell out the results in the case where
the target group $G$ is of \emph{tube type}; this is justified by Theorem
\ref{ZariskiClosure}. We will also assume that $G$ is \emph{adjoint simple}.\\

We now fix an adjoint simple Hermitian Lie group $G$ of tube type and denote by
$\widehat{G} =
\widetilde{G}/\pi_1(G)^{tor}$ the unique central $\ZZ$-extension of $G$. Then
causal geometry gives rise to a bi-invariant partial order on $\widehat{G}$
(see 
\cite{BenSimon_Hartnick_JOLT} for a discussion of this and various related
bi-invariant partial orders on Lie groups). A prototypical example arises from
the action of $G = \PU(1,1)$ on the boundary of the 
Poincar\'e disk $\DD$; this action lifts to an action of the universal covering
$\widehat{G} = \widetilde{\PU(1,1)}$ on $\RR$, hence induces a bi-invariant
partial order on $\widehat{G}$ by setting
\[g \leq h :\Leftrightarrow \forall x\in \RR:\, g.x \leq h.x.\]
In the general case one utilizes the fact that by the tube type assumption there
exists a unique pair $\pm \mathcal C$ of $G$-invariant causal structures on
the Shilov boundary $\check S$ of the bounded symmetric domain associated with
$G$ (see \cite{Kaneyuki}). Here, by a causal structure $\mathcal C$ we mean a
family of \emph{closed} cones $\mathcal C_x \subset T_x \check S$, and
invariance is understood in the sense that $g_*\mathcal C_x = \mathcal C_{gx}$.
The causal structures $\pm \mathcal C$ lift to $\widehat{G}$-invariant causal
structures on the universal covering $\check R$ of $\check S$, which in turn
induce a pair of mutually inverse (closed) partial orders on $\check R$ via causal curves. Let us denote by 
$\preceq$ the partial order which is compatible with the orientation given by the K\"ahler class. 
We then obtain a bi-invariant partial order on $\widehat{G}$ by
setting
\[g \leq_{\widehat{G}} h :\Leftrightarrow \forall x\in \check R:\, g.x \preceq
h.x.\]
The \emph{dominant set} $\widehat{G}^{++}$ (in the sense of \cite{Eliashberg_Polterovich, BenSimon_Hartnick_Comm}) of this bi-invariant order is given by the formula
\[\widehat{G}^{++} := \{g \in \widehat{G}\,|\,\forall h \in G\exists n \in
\mathbb N: \, g^n \geq_{\widehat{G}} h\},\]
We provide the following simple description in terms of the causal structure:
\begin{thm}
If $\widehat{G}$ is of tube type then
\[\widehat{G}^{++} = \{g \in \widehat{G}\,|\,\forall x \in \check R:\, g.x \succ
x\}.\] 
\end{thm}
We now provide an interpretation of weakly-maximal representations in terms of
dominant sets. 
Let  $\Sigma_{g,n}$ be a compact oriented surface of genus $g$
with $n$ boundary components. We always assume that $\chi(\Sigma)_{g,n} < 0$ so
that there exists a hyperbolization
$\rho:\Gamma_{g,n}:=\pi_1(\Sigma_{g,n})\to\PU(1,1)$. If $n\geq1$, then
$\Gamma_{g,n}$ is a free group, hence $\rho$ admits a lift
$\widetilde\rho:\Gamma_{g,n}\to\widetilde{\PU(1,1)}$
whose restriction to the group of homologically trivial loops
$\Lambda_{g,n}:=[\Gamma_{g,n},\Gamma_{g,n}]$ is unique. In particular, the
translation number quasimorphism on $\widetilde{\PU(1,1)}$ pulls back to a
quasimorphism $f_{\Sigma_{g,n}}$ on $\Lambda_{g,n}$. It turns out that this
quasimorphism is independent of the choice of hyperbolization $\rho$; in fact it
admits a topological description in terms of winding numbers \cite{Huber}. In
the case in which $n=0$, one cannot perform this construction on
$\Gamma_{g, 0}$,
but one has to pass to the central extension $\overline{\Gamma_{g, 0}}$
that corresponds to the generator of $\h^2(\Gamma_{g, 0},\ZZ)$ or,
equivalently, 
can be realized as the fundamental group of the $S^1$-bundles over $\Sigma_g$
of Euler number one. One then obtains in the same way as above a canonical
quasimorphism $f_{\Sigma_{g,0}}$ on 
$\Lambda_{g,0}:=[\overline\Gamma_{g,0},\overline\Gamma_{g, 0}]$.
We emphasize that the quasimorphism
$f_{\Sigma_{g,n}}$ depends on the topological surface $\Sigma_{g,n}$, not just the abstract group
$\Gamma_{g,n}$.
\begin{thm}\label{GeometricMain}
Let $G$ be an adjoint simple Hermitian Lie group of tube type and let
$\widehat{G}$,
$\widehat{G}^{++}$ as above. Let $\Sigma_{g,n}$ be a surface
of negative Euler characteristic and $\Gamma_{g,n}:=\pi_1(\Sigma_{g,n})$. Then a
representation $\rho: \Gamma_{g,n} \to G$
is weakly maximal with $T(\rho) \neq 0$ iff for the unique lift $\widetilde\rho:\Lambda_{g, n}\to\widehat{G}$ there exists $N > 0$ such that
\begin{eqnarray}\label{PositiveCausal}
f_{\Sigma_{g,n}}(\gamma) > N \Rightarrow
\widetilde{\rho}(\gamma) \in
\widehat{G}^{++} \quad (\gamma \in \Lambda_{g, n}).
\end{eqnarray}
\end{thm}
\begin{rem}
If we define a family of partial orders $\leq_N$ on $\Lambda_{g, n}$ by
\[g <_N h :\Leftrightarrow f_{\Sigma_{g,n}}(g^{-1}h) > N,\]
and a partial order on $\widehat{G}$ by
\[g \leq_{++} h :\Leftrightarrow g^{-1}h \in \widehat{G}^{++}\]
then the conclusion can be rephrased by saying that $\widetilde{\rho}$ is
order-preserving with respect to some $\leq_N$ and
$\leq_{++}$.
\end{rem}

\section{Comparison to other classes of representations}
Turning back to the general theme of studying subsets of the representation
variety we describe basic properties of the set 
$\hom_{wm}(\pi_1(\Sigma), G)$ of weakly maximal representations, and relate it
to other geometrically meaningful subsets of $\hom(\pi_1(\Sigma), G)$. We will
denote by 
$\hom^*_{wm}(\pi_1(\Sigma), G) \subset \hom(\pi_1(\Sigma), G)$ the subset of
weakly maximal representations with nonzero Toledo number. Also we denote by
$\hom_{di}(\pi_1(\Sigma), G)$ the set of discrete and faithful representations.
By Theorem~\ref{thm_intro: discrete_faithful} we have a chain of inclusions 
\begin{equation}\label{Inclusions1}
\hom_{max}(\pi_1(\Sigma), G) \subset \hom^*_{wm}(\pi_1(\Sigma), G) \subset 
\hom_{di}(\pi_1(\Sigma), G).
\end{equation}
The sets on the right \cite{Goldman_Millson} and on the left are closed; if
$\partial \Sigma = \emptyset$ the left one is also open
\cite{Burger_Iozzi_Wienhard_toledo}. We are able to show:
\begin{thm}\label{thm_intro:closed}
The set $\hom_{wm}(\pi_1(\Sigma), G) \subset \hom(\pi_1(\Sigma), G)$ is closed.
\end{thm}
Combining this with 
Theorem~\ref{thm_intro:integer} we then obtain:
\begin{cor}
The set $\hom^*_{wm}(\pi_1(\Sigma), G) \subset \hom(\pi_1(\Sigma), G)$ is
closed.
\end{cor}
Thus \eqref{Inclusions1} is a chain of \emph{closed} subsets of the
representation variety. In the case where $\Sigma$ is a closed
surface we can refine this chain further: It has been established in
\cite{Burger_Iozzi_Labourie_Wienhard, Burger_Iozzi_Wienhard_anosov} that maximal
representations are Shilov-Anosov in the sense of
\cite{Guichard_Wienhard_anosov, Labourie_anosov}. 
Concerning the (open) set 
$\hom_{\cs-An}(\pi_1(\Sigma), G))$ of all Shilov-Anosov representations
we establish the following:
\begin{thm}
Assume that $\partial \Sigma = \emptyset$ and that $G$ is a Lie group of tube
type. Then 
\begin{equation}
\overline{\hom_{\cs-An}(\pi_1(\Sigma), G)}  \subset \hom_{wm}(\pi_1(\Sigma), G).
\end{equation}
\end{thm}
For a \emph{closed} surface $\Sigma$ and a Hermitian group $G$ of tube type we
thus end up with the following diagram. Here we denote by
$\hom^*_{\epsilon}(\pi_1(\Sigma), G))$ the set of representations in 
$\hom_{\epsilon}(\pi_1(\Sigma), G))$ of nonzero Toledo number. We also denote by
$\hom_{red}(\pi_1(\Sigma), G))$ the set of representations with reductive
Zariski closure and by
$\hom_{Hitchin}(\pi_1(\Sigma), G))$ the Hitchin component in case $G$ is locally
isomorphic to 
$\Sp(2n,\RR)$ (and the empty set otherwise). Then we have the following
inclusions:
\[\begin{xy}\xymatrix{
                      &                 & \hom_{\cs-An} \ar@{}[r]|-*[@]{\subset}&
\hom_{wm}&\\
\hom_{Hitchin}  \ar@{}[r]|-*[@]{\subset} &\hom_{max} \ar@{}[r]|-*[@]{\subset} \ar@{}[d]|-*[@]{\cap}&
\hom^*_{\cs-An} \ar@{}[r]|-*[@]{\subset} \ar@{}[u]|-*[@]{\cup}& \hom^*_{wm} \ar@{}[r]|-*[@]{\subset} \ar@{}[u]|-*[@]{\cup}&
\hom_{di}\\
& \hom_{tight} \ar@{}[r]|-*[@]{\subset} & \hom_{red}
}\end{xy}\]
If $\Sigma$ is allowed to have boundary, then the relations between the various
subsets of the representation variety is more complicated.
\def\cprime{$'$} \def\cprime{$'$}

\end{document}